\documentclass[12pt]{article}
\title{ The Logarithmic Fib-Binomial Formula }
\author{A.K.Kwa\'sniewski\\  
\\ Higher School of Mathematics and Applied Informatics\\
PL-15-021 Bia{\l}ystok, ul. Kamienna 17, POLAND
\\e-mail: kwandr@uwb.edu.pl}

\usepackage{amsmath,amsthm}
\usepackage{graphicx}
\chardef\bslash=`\\ 
\newtheorem{defn}{Definition}[section]
\newtheorem{obs}{Observation}[section]
\newtheorem{nt}{Notation}[section]

\newtheorem{prop}{Proposition}[section]

\begin{document}
\maketitle

\begin{abstract}
Steven Romans Logarithmic Binomial Formula analogue has been found
and is presented here also for the case of  fibonomial
coefficients - which recently have been given a combinatorial
interpretation by the present author.
\end{abstract}

\textit{to appear in} volume 9 (1),\textit{Advanced Studies in
Contemporary Mathematics}

\section{Introduction}
The aim of this note is to find out  - as in \cite{14,15} of
Steven Roman - the form of "Fib-corresponding" Logarithmic
Fib-Binomial Formula. In \cite{14,15} Steven Roman introduced The
Logarithmic Binomial Formula :
$$\lambda_{n}^{(t)}(x+a)=\sum_{k\geq 0}\left[ {n \atop k}\right]\lambda^{(t)}_{n-k}(a)x^{k};
\;\;\;t=0,1;\;\;\;|x|<a;\;\;\;n \in {\bf Z}$$
where
$$[n]=\left\{{n\;\;\;\;n\neq 0 \atop 1\;\;\;\;n=0}\right.$$
 and the  Roman factorial is given by
 \begin{equation}\label{1.1}
 [n]!=\left\{{n!\;\;\;\;\;\;\;\;\;\;n\geq 0 \atop
 \frac{(-1)^{n+1}}{(-n-1)!}\;\;\;n<0}\right.
 \end{equation}
while hybrid binomial coefficients \cite{10}  (Roman coefficients
\cite{11}) read:
\begin{equation}\label{1.2}
\left[{n \atop k}\right]=\frac{[n]!}{[k]![n-k]!}
\end{equation}
 One may show that  (Propositions 3.2 , 4.1, 4,2,
4.3  in \cite{11} )
\begin{equation}\label{1.3}
\left[{n \atop k}\right]=\left\{ \begin{array}{lll}
\binom{n}{k}&&n,k\geq 0\\
\\
(-1)^{k}\binom{-n-1+k}{k}&&k\geq 0 \geq n\\
\\
(-1)^{k+n}\binom{-k-1}{n-k}&&0>n\geq k\\
\\
(-1)^{(n+k)}\left[\Delta^{n}\frac{1}{x-k}\right]_{k=0}&&k>n\geq 0
\end{array}\right.
\end{equation}
\textrm{}\\
$$\left[{0 \atop k}\right]=\left[ {0 \atop -k}\right]=\frac{(-1)^{k+1}}{k}$$
\textrm{}\\
\begin{equation}\label{1.4}
\left[{n \atop k}\right]=\left[ {n \atop
n-k}\right],\;\;\;\;\;\;\; \left[{n \atop j}\right]\cdot \left[ {j
\atop k}\right]=\left[{n \atop k}\right]\cdot \left[ {n-j \atop
j-k}\right]
\end{equation}
\textrm{}\\
\begin{equation}\label{1.5}
\left[{n \atop k}\right]=\left[ {n-1 \atop k-1}\right]+\left[{n-1
\atop k}\right]
\end{equation}
\textrm{}\\
 As seen from the above  the hybrid binomial
coefficients (Roman coefficients) are the intrinsic natural
extension  of binomial coefficients .

 The Logarithmic Binomial Formula  extends the notion of binomiality of polynomials as  used
in the Generalized Umbral Calculus  (see Chapter 6 in \cite{16}
for functional formulation and see \cite{8} for  abundant
references on Finite Operator Calculus of Rota formulation).

The great invention of Steven Roman - among others - relies on the
fact that the real i.e. $R$-linear span $L$  of  the basis
functions (harmonic logarithms-see: Proposition 4.1 in \cite{14} )
$$L=span\left| \left\{ \lambda_{n}^{(t)}\right\}_{n\in {\bf Z},\;t=0,1}\right|$$
allows the Fundamental Theorem of Calculus to hold on  $L$ i.e.
$D^{-1}D = DD^{-1} = id_{L}$. Here $D^{-1}$  depending on whether
$t=0$ or $t=1$ acts as follows:
$$D^{-1}=\int_{0}^{x}\;\;\;on\;\;\;\lambda_{n}^{(0)}\;\;\;n\neq -1,
\;\;n\in {\bf Z}\;\;\;and\;\; gives\;\;0\;\;for\;\;n=-1$$
\textrm{}\\
$$and\;\;\;D^{-1}=\int_{1}^{x}\;\;\;\;on\;\;\;\lambda_{n}^{(1)};\;\;\;n\in {\bf Z}$$

\section{Fibonomial Coefficients}
In \cite{5}  Fibonomial coefficients \cite{12,2,3,4} have been
given a combinatorial interpretation as counting the number of
finite "birth-selfsimilar" subposets of an infinite poset. We
shall  use here the following notation: {\em Fibonomial
coefficients} are defined as
$\binom{n}{k}_{F}=\frac{F_{n}!}{F_{k}!F_{n-k}!}$ or  - usefully
for our purpose here $\binom{n}{k}_{F}\equiv
\frac{n^{\underline{k}}_{F}}{k_{F}!}$ where we make an analogy
driven \cite{8} identifications: ($n>0$), $n_{F}\equiv F_{n}\neq
0$,$n_{F}!\equiv n_{F}(n-1)_{F}(n-2)_{F}(n-3)_{F}\ldots
2_{F}1_{F};\;\;0_{F}!=1$;\\
$n_{F}^{\underline{k}}=n_{F}(n-1)_{F}\ldots (n-k+1)_{F}$. This is
the appropriate specification of notation from \cite{8} for the
purpose Fibonomial Finite Operator Calculus case investigation
(see Example 2.1 in \cite{9}).

Let us now introduce an infinite poset $P$ (for further details
see: \cite{5}) via its finite part subposet $P_{m}$ Hasse diagram
to be continued ad infinitum in an obvious way as seen from the
figure below.  It looks like the Fibonacci tree with a specific
"cobweb": see Figure 1. One sees that the  $P_{m}$  is the
subposet of $P$ consisting of points up to $m$-th level points
$$\bigcup_{s=1}^{m}\Phi_{s}\;\; ;\;\Phi_{s} \;\;is\;\; the\;
\;set\;\;of\;\;elements\;\;of\;\;the\;\;s-th\;\;level$$

How many $P_{m}$`s  \underline{rooted at the $k$-th} level might
be found ?

 We answer this question in the following sequence of
observations right after Figure 1.

\begin{center}

\includegraphics[width=75mm]{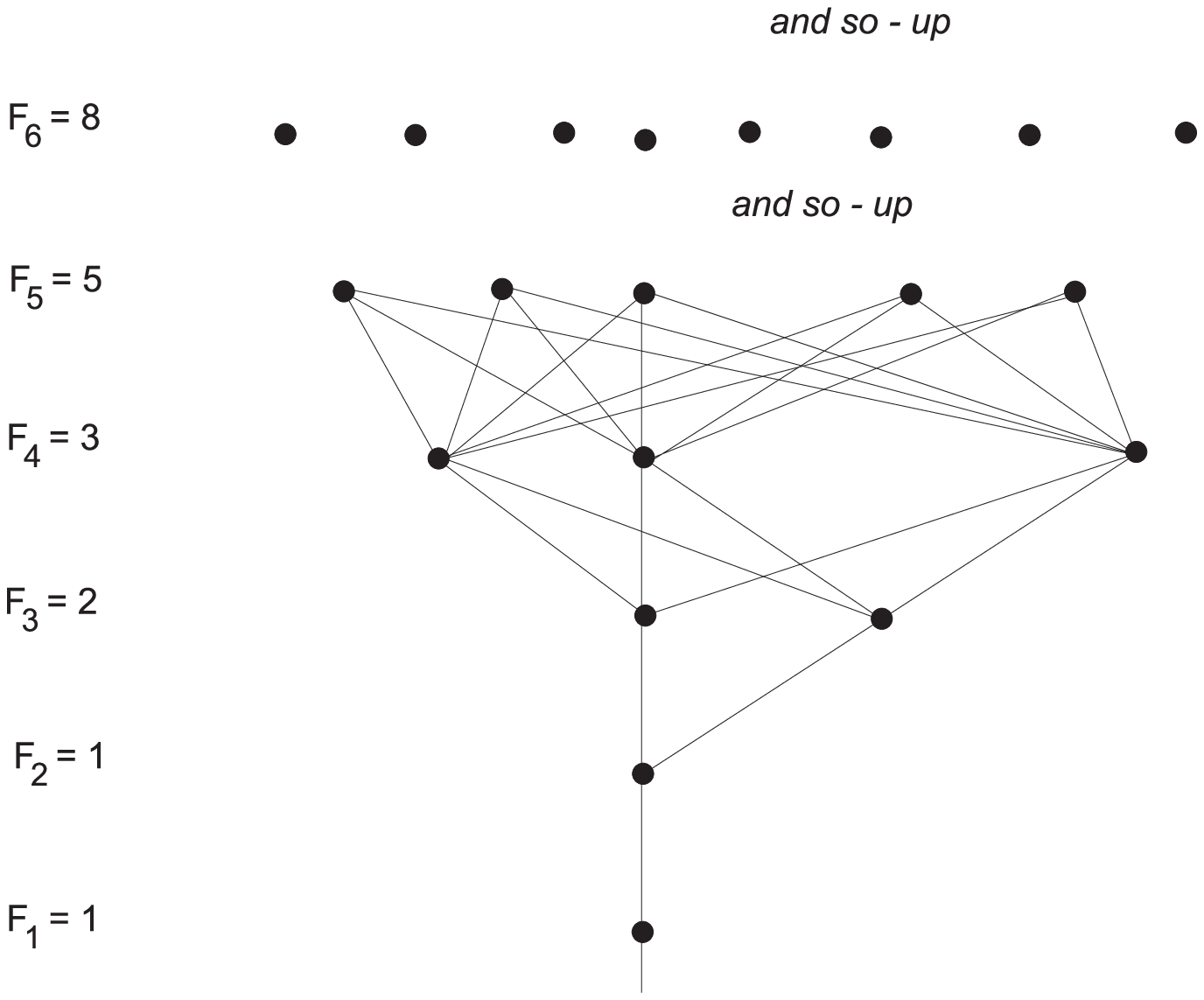}

\vspace{2mm}

\noindent {\small Fig.~1. The construction of the  Fibonacci
"cobweb"  poset} \end{center}
 \begin{obs}
 The number of maximal chains starting from the root  (level  $F_{1}$ )
  to reach any point at the $n$-th level  labeled by  $F_{n}$  is equal
  to $n_{F}!$
  \end{obs}
 \begin{obs}
 The number of maximal chains starting from the level labeled by  $F_{k}$ to
reach any point at the $n$-th level labeled by  $F_{n}$  is equal
to $n_{F}^{\underline{m}},\;\;(n = k+m )$.
\end{obs}
\begin{obs}
Let $n = k+m$. The number of subposets $P_{m}$ rooted at the level
labeled by $F_{k}$ and ending at the $n$-th level labeled by
$F_{n}$ is equal to
$$ \binom{n}{m}_{F}=\binom{n}{k}_{F}=\frac{n^{\underline{k}}_{F}}{k_{F}!}.$$
\end{obs}

\section{The Logarithmic Fibonomial Case}
We shall now adopt the $\ast_{\psi}$  product  formalism \cite{9}
(see also Appendix  in  \cite{6}  and \cite{7}) to the Fibonomial
case with
$\exp_{F}\left\{x\right\}=\sum_{k=0}^{\infty}\frac{x^{k}}{k_{F}!}$
defining the $F$-exponential series.

\subsection{$\ast_{F}$ product}
Let   $n >0$ and  let $\partial_{F}$ be a linear operator acting
on formal series  and defined accordingly by
$\partial_{F}x^{n}=n_{F}x^{n-1};\;\;n\geq
0,\;\;\partial_{F}x^{0}=0$ .

 We shall call the $F$-multiplication the new  $\ast_{F}$
 product of functions or formal series  specified below.
 \begin{nt}\label{not1}
 $x\ast_{F}x^{n}=\hat{x}_{F}(x^{n})=\frac{(n+1)}{(n+1)_{F}}x^{n+1};\;\;n\geq
 0$ hence $x\ast_{F}1=x$ and $x\ast_{F}\alpha 1=\alpha 1\ast_{F}x=
 x\ast_{F}\alpha =\alpha \ast_{F}x=\alpha x;\;\;\;\forall x,\alpha \in {\bf
 R},\;\;f(x)\ast_{F}x^{n}=f(\hat{x}_{F})x^{n}$.

 For $k\neq n$ $x^{n}\ast_{F}x^{k}\neq x^{k}\ast_{F}x^{n}$ as well
 as $x^{n}\ast_{F}x^{k}\neq x^{n+k}$ - in general.
 \end{nt}
 \begin{defn}
  With  Notation \ref{not1} adopted   define  the  $\ast_{F}$
powers of  $x$  according to
$$x^{n^{\ast}_{F}}\equiv x\ast_{F}x^{(n-1)^{\ast}_{F}}=
\hat{x}\left(
x^{(n-1)^{\ast}_{F}}\right)=x\ast_{F}x\ast_{F}\ldots\ast_{F}x=\frac{n!}{n_{F}!}x^{n};\;\;\;n\geq
0.$$
\end{defn}
Note that
$x^{n^{\ast}_{F}}\ast_{F}x^{k^{\ast}_{F}}=\frac{n!}{n_{F}!}x^{(n+k)^{\ast}_{F}}
\neq
x^{k^{\ast}_{F}}\ast_{F}x^{n^{\ast}_{F}}=\frac{k!}{k_{F}!}x^{(n+k)^{\ast}_{F}}$
for $k\neq n$ and $x^{0^{\ast}_{F}}=1$.

This noncommutative $F$-product $\ast_{F}$ is devised so as to
ensure the observations below.
\begin{obs}\label{obs3.1} \textrm{}\\
\renewcommand{\labelenumi}{(\alph{enumi})}
\begin{enumerate}
\item $\partial_{F}x^{n^{\ast}_{F}}=nx^{(n-1)^{\ast}_{F}};\;\;
n\geq0$;
 \item $\exp_{F}[\alpha x]\equiv \exp\left\{
\alpha\hat{x}_{F}\right\}{\bf 1};$
 \item $\exp[\alpha
x]\ast_{F}\left\{\exp_{F}\left\{\beta\hat{x}_{F}\right\}1\right\}=
\exp_{F}\left\{[\alpha+\beta]\hat{x}_{F}\right\}{\bf 1}$;
 \item $\partial_{F}(x^{k}\ast_{F}x^{n^{\ast}_{F}})=
(Dx^{k})\ast_{F}x^{n^{\ast}_{F}}+x^{k}\ast_{F}(\partial_{F}x^{n^{\ast}_{F}})$;
\item Leibniz rule\;\;\;\;\;
$\partial_{F}(f\ast_{F}g)=(Df)\ast_{F}g+f\ast_{F}(\partial_{F}g)$;
$f,g$- formal series; \item
$f(\hat{x}_{F})g(\hat{x}_{F})1=f(x)\ast_{F}\tilde{g};\;\;
\tilde{g}(x)=g(\hat{x}_{F}){\bf 1}$.
\end{enumerate}
\end{obs}
\subsection{$F$-Integration}
 Let :  $\partial_{0}x^{n}=x^{n-1}$. The linear operator $\partial_{0}$ is
 identical with  divided difference operator. Let $\hat{Q}f(x)=f(qx)$.
 Recall  that to the Jackson $\partial_{q}$  derivative \cite{8}  there corresponds
 the $q$-integration  which is a right inverse operation to "$q$-difference-ization".
 Namely \cite{8}
\begin{equation}\label{3.2}
F(z):=\left(
\int_{q}\varphi\right)(z):=(1-q)z\sum_{k=0}^{\infty}\varphi\left(
q^{k}z\right)q^{k}
\end{equation}
\begin{equation}\label{3.3}
F(z)\equiv\left(\int_{q}\varphi\right)(z)=(1-q)z\left(
\sum_{k=0}^{\infty}q^{k}\hat{Q}^{k}\varphi\right)(z)=\left((1-q)\hat{z}\frac{1}{1-q\hat{Q}}\varphi\right)(z)
\end{equation}
where $(\hat{z}\varphi)(z)=z\varphi(z)$.

 Of course
\begin{equation}\label{3.4}
\partial_{q}\circ \int_{q}=id
\end{equation}
  as
\begin{equation}\label{3.5}
\frac{1-q\hat{Q}}{1-q}\partial_{0}\left(
(1-q)\hat{z}\frac{1}{1-q\hat{Q}}\right)=id
\end{equation}
Naturally (\ref{3.5}) might serve to define a right inverse to
Jackson's "$q$-difference-ization"
$\left(\partial_{q}\varphi\right)(x)=\frac{1-q\hat{Q}}{1-q}\partial_{0}\varphi(x)$
and consequently the "$q$-integration " as represented by
(\ref{3.2}) and (\ref{3.3}. As it is well known  the definite
$q$-integral is an numerical approximation of  the definite
integral obtained in the $q\longrightarrow 1$ limit.

Finally we introduce the analogous representation  for
$\partial_{F}$ difference-ization
\begin{equation}\label{3.6}
\partial_{F}=\hat{n}_{F}\partial_{0};\;\;\;\hat{n}_{F}x^{n-1}=n_{F}x^{n-1};\;\;\;n\geq
1.
\end{equation}
Then
\begin{equation}\label{3.7}
\int_{F}x^{n}=\left(
\hat{x}\frac{1}{\hat{n}_{F}}\right)x^{n}=\frac{1}{(n+1)_{F}}x^{n+1};\;\;\;n\geq
0
\end{equation}
 and of course
\begin{equation}\label{3.8}
\partial_{F}\circ \int_{F}=id.
\end{equation}
  Naturally  ($\int_{F}\equiv \int d_{F}$)
  $$\partial_{F}\int_{a}^{x}f(t)d_{F}t=f(x).$$
The formula of   "per partes" $F$-integration  is easily
obtainable from  Observation (\ref{obs3.1}) and it  reads:
\begin{equation}\label{3.9}
\int_{a}^{x}(f\ast_{F}\partial_{F}g)(t)d_{F}t=\left[
(f\ast_{F}g)(t)\right]^{x}_{a}-\int_{a}^{x}(Df\ast_{F}g)(t)d_{F}t.
\end{equation}

Now in order to have  $\partial_{F}^{-1}$- an  $F$-analogue  of
$D^{-1}$ as in \cite{13,14}( thus causing the fundamental Theorem
of Calculus to hold for $\partial_{F}$-difference-ization and
$\int_{F}$- integration on some linear space $L_{F}$ being the
linear span of "$F$-{\em harmonic logarithms}" )  - we shall
proceed exactly as Steven Roman  in \cite{14,15}.

\section{The Logarithmic Fib-Binomial Formula}
As in  \cite{14,15} of  Roman  -  we have also  {\em The
Logarithmic Fib-Binomial Formula} (see: Propositions
\ref{prop4.1}, \ref{prop4.2} below):
$$\phi_{n}^{(t)}(x+_{F}a)\equiv \left[ \exp\left\{ a\partial_{F}\right\} \phi_{n}^{(t)}\right](x)
=\sum_{k \geq 0}\left[{n \atop
k}\right]_{F}\phi_{n-k}^{(t)}(a)x^{k}\;\;\;t=0,1;\;|x|<a;\;n\in
{\bf Z}$$
 where ( more on "$+_{F}$ " see \cite{8,9})
 $$[n_{F}]=\left\{\begin{array}{lll}
 n_{F}&&n\neq 0\\
 1&&n=0\end{array}\right.$$
 and the  Roman Fib-factorial is given by
\begin{equation}\label{4.1}
[n_{F}]!=\left\{\begin{array}{lll} n_{F}!&&n\geq 0\\
\frac{(-1)^{n+1}}{(-n-1)_{F}!}&&n<0\end{array}\right.
\end{equation}
while {\em Fib-hybrid binomial coefficients} or Roman
Fib-coefficients (see: \cite{10,11}) read:
\begin{equation}\label{4.2}
\left[ {n \atop k}\right]_{F}=\frac{[n]!}{[k]![n-k]!}
\end{equation}
One observes   (as in  Propositions 3.2 , 4.1, 4,2, 4.3  in
\cite{11} ) that:
\begin{equation}\label{4.3}
\left[ {n \atop k}\right]_{F}=\left\{\begin{array}{lll}
\binom{n}{k}_{F}&&n,k\geq 0\\
\\
(-1)^{k}\binom{-n-1+k}{k}_{F}&&k\geq 0>n\\
\\
(-1)^{k+n}\binom{-k-1}{n-k}_{F}&&0>n\geq k\\
\\
(-1)^{(n+k)}\left[\Delta_{F}^{n}\frac{1}{x-k}\right]_{k=0}&&k>n\geq0
\end{array}\right.
\end{equation}
\textrm{}\\
$$\left[{0 \atop k}\right]_{F}=\left[ {0 \atop -k}\right]_{F}=\frac{(-1)^{k+1}}{k_{F}}$$
\textrm{}\\
\begin{equation}\label{4.4}
\left[{n \atop k}\right]_{F}=\left[ {n \atop
n-k}\right]_{F},\;\;\;\;\;\;\; \left[{n \atop j}\right]_{F}\cdot
\left[ {j \atop k}\right]_{F}=\left[{n \atop k}\right]_{F}\cdot
\left[ {n-j \atop j-k}\right]_{F}
\end{equation}
\textrm{}\\
\begin{equation}\label{4.5}
\left[{n \atop k}\right]_{F}=\left[ {n-1 \atop
k-1}\right]_{F}+\left[{n-1 \atop k}\right]_{F}
\end{equation}
\textrm{}\\
where (see: pp.333-334 in  \cite{8} )
$$\Delta_{F}=\exp_{F}\left\{ \partial_{F}\right\}-id.$$

Fib-Roman coefficients (as  seen  from the above) are then also
natural "relative" of binomial coefficients among the family of
$\psi$ - binomial ones  \cite{8}  (consult also Example 2.1 in
\cite{9}).

 The Logarithmic Fib-Binomial Formula  extends the
notion of binomiality of \underline{polynomials} as  used in the
Generalized Umbral Calculus (see Chapter 6 in \cite{16} for
functional formulation and see \cite{8} for  abundant  references
on Finite Operator Calculus of Rota formulation)- to sequences of
\underline{functions} - (compare with \cite{1}).

Here the importance of the great invention of Steven Roman - among
others - relies on the fact that the $R$-linear span $L_{F}$ of
now basis {\em Fib-harmonic logarithms} functions
$$\left\{ \phi_{n}^{(t)}\right\}_{n\in {\bf Z},\;t=0,1},\;\;\;\;\;L_{F}
=span\left| \left\{ \phi_{n}^{(t)}\right\}_{n\in {\bf
Z},\;t=0,1}\right|$$ allows the Fundamental Theorem of Calculus to
hold also on $L_{F}$, i.e.\\
$\partial_{F}^{-1}\partial_{F}=id_{L_{F}}$ for $\partial_{F}$-
difference-ization and  $\int_{F}$ - integration acting on \\a
linear space $L_{F}$ being the linear span of "{\em Fib-harmonic
logarithms}". Here anti-difference-ization operator
$\partial_{F}^{-1}$ - depending on whether $t = 0$  or  $t = 1$\\
- acts as follows on  {\em Fib-harmonic logarithm} functions :
$$\partial_{F}^{-1}=\int_{0}^{x}d_{F}\;\;\;on\;\;\;\phi_{n}^{(0)}\;\;\;n\neq -1,
\;\;n\in {\bf Z}\;\;\;and\;\; gives\;\;0\;\;for\;\;n=-1$$
$$and\;\;\;\partial_{F}^{-1}=\int_{1}^{x}d_{F}\;\;\;\;on\;\;\;\phi_{n}^{(1)};\;\;\;n\in {\bf Z}$$
 Let us define these {\em Fib-harmonic logarithms}
 $$\left\{ \phi_{n}^{(t)}\right\}_{n\in {\bf Z},\;t=0,1}$$
-(see  Proposition 2.2 in \cite{14}) - as solutions of {\em
Fib-harmonic} $t$-binomiality conditions. Thus  {\em Fib-harmonic
logarithm \underline{functions}} are unique solutions of {\em
Fib-harmonic} $t$-binomiality conditions; $t = 0,1$  (\ref{4.6}) -
(compare with  \cite{1} and relaxation Lemma 2.12 therein):
\begin{equation}\label{4.6}
\begin{array}{l}
1)\;\;
\phi_{0}^{(0)}(x)=1,\quad\quad\;\;\;2)\;\;\phi_{n}^{(0)}(0)=0,\;\;n\ni
{\bf Z}\backslash\{ 0\},\\
\\3)\;\;\partial_{F}\phi_{n}^{(0)}=[n_{F}]\phi_{n-1}^{(0)},\;\;\;n\in
{\bf Z}\\ \\ \\
1)\;\;\phi_{0}^{(1)}(x)=\ln x,\quad
\quad2)\;\;\phi_{n}^{(1)}(x)\;has\; no\;constant\;term,\;\;n\in
{\bf Z},\\
\\3)\;\;\partial_{F}\phi_{n}^{(1)}=[n_{F}]\phi_{n-1}^{(1)},\;\;\;n\in
{\bf Z}
\end{array}
\end{equation}
The {\em Fib-harmonic} $t$-binomiality conditions;  $t = 0,1$
(\ref{4.6}) yield \cite{14}  what follows:
\begin{prop}\label{prop4.1}
\begin{multline*}
\phi_{n}^{(0)}(x)=\left\{\begin{array}{lll}
x^{n}&&n\geq 0\\
0&&n<0\end{array}\right.,\;\;\;
\phi_{n}^{(1)}(x)=\left\{\begin{array}{lll} x^{n}(\ln
x-f_{n})&&n\geq 0\\
x^{n}&&n<0\end{array}\right. , \\
\;\;f_{0}=0,\;\;f_{n}=1+\frac{1}{2_{F}}+\frac{1}{3_{F}}+\ldots+\frac{1}{n_{F}},\;\;n\in
{\bf N}
\end{multline*}
\end{prop}
We shall call $f_{n},\;n\in {\bf N}$ the {\em
\underline{Fib-harmonic numbers}} ($f_{0}=0$), (see: \cite{13}).
\begin{prop}\label{prop4.2}
The linear  anti-difference-ization unique operator\\
$\partial_{F}^{-1}:L_{F}\longrightarrow L_{F}$;
$\partial_{F}^{-1}\partial_{F}=id_{L_{F}}$ is given by
\textrm{}\\
$$\partial_{F}^{-1}\phi_{n}^{(0)}=\left\{\begin{array}{lll}
\frac{1}{[n+1]_{F}}\phi_{n+1}^{(0)}&&n\neq -1\\
\\
0&&n=-1\end{array}\right.,\;\;\;\partial_{F}^{-1}\phi_{n}^{(1)}=\frac{1}{[n+1]_{F}}\phi_{n}^{(1)},\;\;n\in
{\bf Z}.$$
\textrm{}\\
\end{prop}
{\bf REMARK.} Instead of  Roman Fib-coefficients and  Roman
Fib-factorial  one may - (replace $F$ by $\psi$)- start to
consider Roman  $\psi$ -coefficients, $\psi$-{\em harmonic
logarithms} etc. However these seemingly might lack any
"reasonable" combinatorial interpretation.

 As the generally useful reading - also for this purpose one recommends here:\\
$[$LR$]$ Loeb D.,  Rota G-C. {\em "Recent Contributions to the
calculus of Finite Differences: a Survay"}  Lecture Notes in Pure
and Appl. Math. vol. \underline{132}(1991), pp. 239-276,  ArXiv:
math.co/ 9502210 V1 9 Feb 1995, see also:\\
http://arxiv.org/list/math.CO/9502

A.M.S. Classification numbers  11B39 , 11B65

 \end{document}